\documentclass[10pt, a4paper]{scrartcl}

\usepackage{latexsym}
\usepackage{amssymb}
\usepackage[nonamelimits]{amsmath}
\usepackage[authoryear, round]{natbib}

\usepackage{color}
\definecolor{myred}{rgb}{0.5,0,0}
\definecolor{myblue}{rgb}{0,0,0.75}
\definecolor{mygreen}{rgb}{0,0.5,0}

\usepackage{ifpdf}

\ifpdf
   \usepackage[pdftex]{graphicx}
   \usepackage[pdftex,
               pdftitle={The Law of Total Odds},
               pdfsubject={},
               pdfkeywords={Total probability, likelihood ratio, Bayes' formula, binary classification, relative odds, 
unbiased estimator, supervised learning, dataset shift},
               pdfauthor={Dirk Tasche},
               pdfstartview=FitR,
               breaklinks=true,
               colorlinks=true,
               citecolor=myred,
               linkcolor=myblue,
               urlcolor=mygreen]{hyperref}
\else
   \usepackage[dvips]{graphicx}
   \usepackage{hyperref}
   \usepackage{rotating}
\fi

\newtheorem{theorem}{Theorem}[section]
\newtheorem{lemma}[theorem]{Lemma}

\newtheorem{proposition}[theorem]{Proposition}

\newtheorem{corollary}[theorem]{Corollary}
\newtheorem{assumption}[theorem]{Assumption}

\newcommand{\qed}{$\Box$}

\numberwithin{equation}{section}

\title{The Law of Total Odds}

\author{%
Dirk Tasche\thanks{E-mail: dirk.tasche@gmx.net\newline
The author currently works at the Prudential Regulation Authority (a
division of the Bank of England). He is also a visiting professor at
Imperial College, London. 
The opinions expressed in this paper are those of the author 
and do not necessarily reflect views of the Bank of England.}}

\date{First version: December 3, 2013\\
This version: February 14, 2014}

\begin{document}

\maketitle

\begin{abstract}
The law of total probability may be deployed in binary classification exercises
to estimate the unconditional class probabilities if the class proportions in the training set 
are not representative of the population class proportions. We argue that this is not
a conceptually sound approach and suggest an alternative based on the 
new law of total odds. We quantify the bias of the total probability estimator
of the unconditional class probabilities and show that the total odds estimator is unbiased. 
The sample version of the total odds estimator is shown to coincide with a maximum-likelihood
estimator known from the literature.
The law of total odds can also be used for transforming
the conditional class probabilities if independent estimates of the unconditional 
class probabilities of the population are available.\\
\textsc{Keywords:} Total probability, likelihood ratio, Bayes' formula, binary classification, relative odds, 
unbiased estimator, supervised learning, dataset shift. 
\end{abstract}

\section{Introduction}
\label{se:intro}

The law of total probability is one of the fundamental building blocks of probability theory. Its elementary 
version states that for an event $A$ and a partition $H_i$, $i\in\mathbb{N}$ of the whole space 
the probability of $A$ can be calculated as
\begin{align}
	\mathrm{P}[A] &= \sum_{i=1}^\infty \mathrm{P}[H_i]\,\mathrm{P}[A\,|\,H_i], \label{eq:basic}
\intertext{where the conditional probabilities $\mathrm{P}[A\,|\,H_i]$ are defined as}
\mathrm{P}[A\,|\,H_i] &= 
	\begin{cases}
\frac{\mathrm{P}[A\cap H_i]}{\mathrm{P}[H_i]}, & \text{if}\ \mathrm{P}[H_i] > 0,\\
0, & \text{if}\ \mathrm{P}[H_i] = 0.
\end{cases}\notag
\end{align}
\citet{kolmogorov1956foundations} calls Eq.~\eqref{eq:basic} the theorem of total probability. It is also
called rule or formula of total probability. Virtually all text books on probability theory mention
Eq.~\eqref{eq:basic} but many authors \citep[e.g.][Chapter~V, Eq.~(1.8)]{feller1968probability} do not name it. 

\citet{feller1968probability} comments on Eq.~\eqref{eq:basic} with the words ``This formula is useful because
an evaluation of the conditional probabilities $\mathrm{P}[A\,|\,H_i]$ is frequently easier than a direct calculation
of $\mathrm{P}[A]$.'' Sometimes it may even be impossible to directly calculate $\mathrm{P}[A]$. In particular, 
this is the case when $\mathrm{P}[A]$ is assumed to be forecast but past observations of occurrences of event $A$ cannot
be relied on because the value of $\mathrm{P}[A]$ might have changed.

Such a situation is likely to be incurred in binary classification exercises where the unconditional (or prior) 
class probabilities
in the training dataset may differ from the class probabilities of the population to which the classifier is applied
\citep[see][for a recent survey of data shift issues in classification]{MorenoTorres2012521}. Typically, 
a classifier produces class probabilities (i.e.\ probabilities of tested examples
to be of -- say -- class $A$) conditional on already known features $H_n$ of the examples. If the unconditional 
distribution of the $H_n$ (i.e.\ the probabilities $\mathrm{P}[H_i]$) is also known, Eq.~\eqref{eq:basic} then can 
be used to make a forecast (or point estimate) of $\mathrm{P}[A]$.

It can be argued, however, that the forecasts of unconditional class probabilities produced this way are biased
(see Section~2.2 of \citealp{Xue:2009:QSC:1557019.1557117}, or \citealp{Tasche2013}, and Proposition~\ref{pr:biased} below). 
This is a consequence of the fact that fundamentally the conditional class probabilities 
$\mathrm{P}[A\,|\,H_i]$ are determined by means of Bayes' formula (assuming $\mathrm{P}[H_i]>0$ and $\mathrm{P}[A]>0$):
\begin{equation}\label{eq:bayes}
\begin{split}
	\mathrm{P}[A\,|\,H_i] & = \frac{\mathrm{P}_0[A]\,\mathrm{P}[H_i\,|\,A]}{\mathrm{P}_0[A]\,\mathrm{P}[H_i\,|\,A] +
			\mathrm{P}_0[A^c]\,\mathrm{P}[H_i\,|\,A^c]}\\
			& = \frac{\mathrm{P}_0[A]}{\mathrm{P}_0[A] +
			\mathrm{P}_0[A^c]\,\frac{\mathrm{P}[H_i\,|\,A^c]}{\mathrm{P}[H_i\,|\,A]}},
\end{split}
\end{equation}
where $A^c$ denotes the event complementary to $A$ and $\mathrm{P}_0[A]$ and $\mathrm{P}_0[A^c]$ 
are the unconditional probabilities of class $A$  and $A^c$ respectively in the training dataset.
The conditional probabilities $\mathrm{P}[H_i\,|\,A]$ and $\mathrm{P}[H_i\,|\,A^c]$ reflect the
distributions of the characteristic features on class $A$ and its complementary class respectively.

On the one hand, Eq.~\eqref{eq:bayes} suggests a potentially unintended impact of the training set class probabilities
on the population class estimates. On the other hand, Eq.~\eqref{eq:bayes} also suggests that an estimate of $\mathrm{P}[A]$
based solely on the conditional likelihood ratio $i \mapsto \lambda_i = \frac{\mathrm{P}[H_i\,|\,A^c]}{\mathrm{P}[H_i\,|\,A]}$ would
avoid this issue.

This paper presents in Theorem~\ref{th:main} below a necessary and sufficient criterion for when it is possible to estimate
a population class probability based on the unconditional distribution of the features of the tested examples and the
conditional likelihood ratio.
The likelihood ratio $\lambda_i$ can also be written as
\begin{equation}\label{eq:odds}
	\lambda_i = \frac{\mathrm{P}[A^c\,|\,H_i]}{\mathrm{P}[A\,|\,H_i]}\, \frac{\mathrm{P}_0[A]}{\mathrm{P}_0[A^c]}.
\end{equation}
By Eq.~\eqref{eq:odds}, $\lambda_i$ can alternatively be described as the ratio of the conditional and
unconditional odds of class $A^c$ or the \emph{relative odds} of class $A^c$. This observation suggests that  
Theorem~\ref{th:main} is called \emph{law of total odds} in analogy to the law of total probability Eq.~\eqref{eq:basic}.

It turns out that the prior class probability estimator suggested by Theorem~\ref{th:main} is the two-class special case of
the maximum likelihood estimator discussed by \citet{saerens2002adjusting}.  Equation~\eqref{eq:unique} from Theorem~\ref{th:main}
has recently been studied in the $n$-class case by \citet[][Eq.~(9)]{duPlessis2014110}. The contributions of this paper
(limited to the case of binary classification)
to the existing literature can be described as follows:
\begin{itemize}
	\item It is shown that the total odds estimator not only solves a prior probability shift problem\footnote{%
	See \citet{MorenoTorres2012521} for the definitions of the various datashift problems.}
	but also, at the same time,
	a combined covariate shift and concept shift problem where only the relative odds are the same for training set and population (or test set).
	\item We demonstrate that the maximum likelihood estimator introduced by \citet{saerens2002adjusting} and studied in more 
	detail by \citet{duPlessis2014110} does not always exist.
	\item We show how to determine conditional class distributions in the population or test set.
	\item It becomes clear that -- in the binary case -- the total odds estimates can be computed by simple numerical root-finding. 
	There is no need to deploy the expectation-maximisation or other more advance iterative algorithms as discussed by
	\citet{saerens2002adjusting}, \citet{Xue:2009:QSC:1557019.1557117} or \citet{duPlessis2014110}.
	\item We provide sharp error bounds for the prior class probability estimate when the covariate shift is ignored. This approach
	is called 'total probability' below.
\end{itemize}


\section{Results}
\label{se:results}

It is useful to consider the use of Eq.~\eqref{eq:basic} for estimating the class probability $\mathrm{P}[A]$ in a more general 
setting.

\begin{assumption}\label{as:general}
$(\Omega, \mathcal{A}, \mathrm{P}_0)$ is a probability space. $\mathcal{H}$ is a sub-$\sigma$-field of $\mathcal{A}$, i.e.\
$\mathcal{H} \subset \mathcal{A}$. $\mathrm{P}_1$ is a probability measure on $(\Omega, \mathcal{H})$ that is absolutely continuous
with respect to $\mathrm{P}_0\bigm|\mathcal{H}$, i.e.\ $\mathrm{P}_1 \ll \mathrm{P}_0\bigm|\mathcal{H}$. $\mathrm{E}_i$ denotes
the expectation operator based on $\mathrm{P}_i$.
\end{assumption}
The interpretation of Assumption~\ref{as:general} is as follows: 
\begin{itemize}
	\item $(\Omega, \mathcal{A}, \mathrm{P}_0)$ is a model that has been fit to historical observations (e.g.\ the training set of a classifier).
	\item $\sigma$-field $\mathcal{H}$ represents the scores produced by the model (classifier) while $\sigma$-field 
	$\mathcal{A}$ additionally contains information regarding the classes of the tested examples.
	\item $(\Omega, \mathcal{H}, \mathrm{P}_1)$ is the outcome of an application of the model to a different set of -- possibly
	more up-to-date -- observations. $(\Omega, \mathcal{H}, \mathrm{P}_1)$  could be a representation of the distribution of  
	the scores produced by the classifier.
	\item The general problem is to extend $\mathrm{P}_1$ to $\mathcal{A}$, by using information from $(\Omega, \mathcal{A}, \mathrm{P}_0)$.
	\item More specifically, the problem might only be to obtain an estimate $\mathrm{P}^\ast_1[A]$ for a fixed event (or class) 
	$A \in \mathcal{A}\backslash\mathcal{H}$,
	as described in Section~\ref{se:intro}.
	However, to make sure that the estimate is meaningful it should be based on a valid model -- which would be an extension of 
	$\mathrm{P}_1$ to any $\sigma$-field containing $A$.
	\item $\mathrm{P}_1 \ll \mathrm{P}_0\bigm|\mathcal{H}$ is a technical assumption that has intuitive appeal, however. For prediction
	based on $(\Omega, \mathcal{A}, \mathrm{P}_0)$ would be pointless if there were events that were possible under $\mathrm{P}_1$ but 
	impossible under $\mathrm{P}_0$.
\end{itemize}
The most obvious extension of $\mathrm{P_1}$ to $\mathcal{A}$ is by means of the conditional probabilities $\mathrm{P}_0[A\,|\,\mathcal{H}]$ 
determined under the measure $\mathrm{P_0}$. Formally, the extension is defined by
\begin{equation}\label{eq:general}
	\mathrm{P}^\ast_1[A] = \mathrm{E}_1\bigl[\mathrm{P}_0[A\,|\,\mathcal{H}]\bigr], \ A \in \mathcal{A}.
\end{equation}
We note without proof that under Assumption~\ref{as:general} $\mathrm{P}^\ast_1$ behaves as we might have expected.
\begin{proposition}\label{pr:noproof}
Under Assumption~\ref{as:general} the set function $\mathrm{P}^\ast_1$ defined by \eqref{eq:general} is a probability
measure on $(\Omega, \mathcal{A})$ with $\mathrm{P}^\ast_1\bigm|\mathcal{H} = \mathrm{P}_1$ and
$\mathrm{P}^\ast_1[A\,|\,\mathcal{H}] = \mathrm{P}_0[A\,|\,\mathcal{H}]$. 
\end{proposition}
Eq.~\eqref{eq:basic} represents the special case of 
Eq.~\eqref{eq:general} where $\mathcal{H} = \sigma(H_n:\,n \in \mathbb{N})$ is 
a $\sigma$-field generated by a countable partition of $\Omega$.

The odds-based alternative to Eq.~\eqref{eq:general} requires more effort and works for single events at a time only.
For $M \subset \Omega$ let
$M^c = \Omega \backslash M$ denote the complement of $M$.
\begin{assumption}\label{as:LR}
Assumption~\ref{as:general} holds. An event $A \in \mathcal{A}$ with $0 < \mathrm{P}_0[A] \stackrel{\text{def}}{=} p_0 < 1$ is fixed.  The two 
conditional distributions $H \mapsto \mathrm{P}_0[H\,|\,A]$ and  $H \mapsto \mathrm{P}_0[H\,|\,A^c]$, $H \in \mathcal{H}$
are absolutely continuous with respect to some $\sigma$-finite measure $\mu$ on $(\Omega, \mathcal{H})$. Denote by $f_A$ and
$f_{A^c}$ the $\mu$-densities of  $\mathrm{P}_0[\cdot\,|\,A]$ and $\mathrm{P}_0[\cdot\,|\,A^c]$ respectively. 
Both $f_A$ and $f_{A^c}$ are positive $\mu$-almost everywhere.
\end{assumption}
The assumption of absolute continuity of the conditional distributions is not really a restriction because one can always 
choose $\mu = \mathrm{P}_0\bigm|\mathcal{H}$. Typically, in practical applications $\mathcal{H}$ is a proper sub-$\sigma$-field of $\mathcal{A}$
and generated by a statistic like a score function. 
It is therefore likely to have $\mu =$ Lebesgue measure on $\mathbb{R}^d$ or $\mu =$ some counting measure.
The assumption of positive densities is more restrictive but intuitive because statistical prediction of events that were impossible in the past does 
not make much sense.

The following proposition provides the general version of Eq.~\eqref{eq:bayes}. We omit its well-known proof.

\begin{proposition}\label{pr:cond.prob}
Under Assumption~\ref{as:LR}, define the conditional likelihood ratio $\lambda_0$ by $\lambda_0 = \frac{f_{A^c}}{f_{A}}$. 
Then it holds that
\begin{itemize}
	\item[(i)] $f = p_0\,f_A + (1-p_0)\,f_{A^c}$ is a $\mu$-density of $\mathrm{P}_0\bigm|\mathcal{H}$, and
	\item[(ii)] $\mathrm{P}_0[A\,|\,\mathcal{H}]$ can be represented as
	$\displaystyle\mathrm{P}_0[A\,|\,\mathcal{H}] = \frac{p_0}{p_0 + (1-p_0)\,\lambda_0}.$
\end{itemize}
\end{proposition}
Consider the special case of $\lambda_0 =1$ in Proposition~\ref{pr:cond.prob}. It holds that
\begin{equation}
	\mathrm{P}_0[\lambda_0 =1] =1 \ \iff\ \mu(f_A \not= f_{A^c}) = 0 \ \iff\ \mathcal{H}\ \text{and}\ A\ \text{are independent}.
\end{equation}
This case is not of much interest for classification problems because it means that $\mathcal{H}$ does not carry any
information with regard to $A$ or $A^c$. We will therefore exclude it from the following discussions. But note that
by the absolute continuity requirement of Assumption~\ref{as:general} $\mathrm{P}_0[\lambda_0 =1] =1$ implies 
$\mathrm{P}_1[\lambda_0 =1] =1$ but $\mathrm{P}_0[\lambda_0 =1] < 1$ in general is not sufficient for $\mathrm{P}_1[\lambda_0 =1] < 1$.

With Proposition~\ref{pr:cond.prob}, we are in a position to state the main result of this note.
Denote by $\mathbf{1}_M$ the indicator function of the event $M$, i.e.\ $\mathbf{1}_M(\omega)=1$
for $\omega \in M$ and $\mathbf{1}_M(\omega)=0$ for  $\omega \in M^c$.
\begin{theorem}[Law of total odds]\label{th:main}
Let Assumption~\ref{as:LR} hold and define the likelihood ratio $\lambda_0$ as in Proposition~\ref{pr:cond.prob}.
Suppose that $\mathrm{P}_1[\lambda_0 =1] < 1$.
\begin{itemize}
	\item[(i)] There exists a solution $0 < p_1 <1$ to the equation
\begin{equation}\label{eq:unique}
	1 = \mathrm{E_1}\left[\frac{1}{p_1 + (1-p_1)\,\lambda_0}\right]
\end{equation}
if and only if $\mathrm{E_1}[\lambda_0] > 1$ and $\mathrm{E_1}[\lambda_0^{-1}] > 1$. If there is a
solution $0 < p_1 <1$ to Eq.~\eqref{eq:unique} it is unique.
\item[(ii)] Let $\mathcal{H}^A = \sigma(\mathcal{H} \cup \{A\})$ denote the $\sigma$-field generated by
$\mathcal{H}$ and $A$. Then it holds that 
$$\mathcal{H}^A = \bigl\{(A\cap H) \cup (A^c\cap G): H, G \in \mathcal{H}\bigr\}.$$
\item[(iii)]  If there is a solution $0 < p_1 <1$ to Eq.~\eqref{eq:unique}
define $\mathrm{P}_1^\ast[B]$ for $B \in \mathcal{H}^A$ by
\begin{equation*}
	\mathrm{P}_1^\ast[B] = \mathrm{E_1}\left[\mathbf{1}_{H}\,\frac{p_1}{p_1 + (1-p_1)\,\lambda_0}\right] +
	\mathrm{E_1}\left[\mathbf{1}_{G}\,\frac{(1-p_1)\,\lambda_0}{p_1 + (1-p_1)\,\lambda_0}\right], 
\end{equation*}
for any representation $(A\cap H) \cup (A^c\cap G)$ of $B$ with $H, G \in \mathcal{H}$.
 Then $\mathrm{P}_1^\ast$ is a probability measure on $\mathcal{H}^A$ with $\mathrm{P}^\ast_1\bigm|\mathcal{H} = \mathrm{P}_1$.
\item[(iv)] The conditional probability $\mathrm{P}_1^\ast[A\,|\,\mathcal{H}]$ is given by
	$$\displaystyle\mathrm{P}_1^\ast[A\,|\,\mathcal{H}] = \frac{p_1}{p_1 + (1-p_1)\,\lambda_0}.$$
\end{itemize}
\end{theorem}

The proof of Theorem~\ref{th:main} is given in Section~\ref{se:proofs} below. Let us
note here instead some observations on Theorem~\ref{th:main}:
\begin{itemize}
	\item The definition of $\mathrm{P}_1^\ast$ and Eq.~\eqref{eq:unique} imply
	$\mathrm{P}_1^\ast[A] = p_1$. Hence we have shown that, by means of Eq.~\eqref{eq:unique}, the
	total odds approach provides a properly modelled population (or test set) estimate of the unconditional
	probability of class $A$ if the condition for likelihood ratio $\lambda_0$ from Theorem~\ref{th:main} (i) is 
	satisfied.
	\item From Proposition~\ref{pr:cond.prob} (ii) and Theorem~\ref{th:main} (iv)
 it follows that
\begin{equation}\label{eq:gen.odds}
	\frac{\mathrm{P}_0[A^c\,|\,\mathcal{H}]}{\mathrm{P}_0[A\,|\,\mathcal{H}]}\,
		\frac{p_0}{(1-p_0)} = \lambda_0 = \frac{\mathrm{P}_1^\ast[A^c\,|\,\mathcal{H}]}{\mathrm{P}_1^\ast[A\,|\,\mathcal{H}]}\,
		\frac{p_1}{(1-p_1)}.
\end{equation}
Hence $\lambda_0$ has an interpretation as relative odds and is the same for both the training set model $\mathrm{P}_0$ and
the population model $\mathrm{P}_1^\ast$. This justifies the naming of Theorem~\ref{th:main}. 
		\item The proof of Theorem~\ref{th:main} (iv) (see Section~\ref{se:proofs}) shows that	
	\begin{equation}
		\mathrm{P}_1^\ast[H\,|\,A] = \mathrm{E_1}\left[\mathbf{1}_{H}\,\frac{1}{p_1 + (1-p_1)\,\lambda_0}\right], \quad H \in \mathcal{H}.
	\end{equation}
	Hence, Eq.~\eqref{eq:unique} ensures that the conditional distribution $H \mapsto \mathrm{P}_1^\ast[H\,|\,A]$ is
	properly normalised.
	\item Violation of the condition for $\lambda_0$ from Theorem~\ref{th:main} (i) could be interpreted
	as evidence that between the observations of $\mathrm{P}_0$ and
	$\mathrm{P}_1$ circumstances  have changed so much that the two models associated with the measures are incompatible. 
	\item In the special case where $\mathcal{H} = \sigma(H_n:\,n \in \mathbb{N})$ is 
a $\sigma$-field generated by a countable partition of $\Omega$, Eq.~\eqref{eq:unique}
reads
\begin{equation}
	1 = \sum_{n=1}^\infty \frac{\mathrm{P}_1[H_n]}{p_1+(1-p_1)\,\frac{\mathrm{P}_0[H_n\,|\,A^c]}{\mathrm{P}_0[H_n\,|\,A]}}.
\end{equation}
 Basically, this is Eq.~(3.11a) of \citet{Tasche2013}, but with a possibly infinite number of `rating grades'.
\end{itemize}
\begin{corollary}\label{co:unique}
The probability measure $\mathrm{P}_1^\ast$ from Theorem~\ref{th:main} is unique in the following sense: If $\widetilde{\mathrm{P}}_1$ 
is any probability measure on $\mathcal{H}^A$ with $\widetilde{\mathrm{P}}_1[A] \in (0,1)$, $\widetilde{\mathrm{P}}_1
\bigm|\mathcal{H} = \mathrm{P}_1$, and
\begin{equation}\label{eq:equiv}
	\frac{\widetilde{\mathrm{P}}_1[A^c\,|\,\mathcal{H}]}{\widetilde{\mathrm{P}}_1[A\,|\,\mathcal{H}]} \,
	\frac{\widetilde{\mathrm{P}}_1[A]}{\widetilde{\mathrm{P}}_1[A^c]}\ = \ \lambda_0,
\end{equation}
then it follows that $\widetilde{\mathrm{P}}_1 = \mathrm{P}_1^\ast$.
\end{corollary}
By Corollary~\ref{co:unique}, a probability measure on $\mathcal{H}^A$  is uniquely determined by 
the marginal distribution on $\mathcal{H}$ and the relative odds with respect to the event $A$.
See Section~\ref{se:proofs} for a proof of the corollary.
 
The real-world estimation exercise from \citet[][Section~4.4]{Tasche2013} shows that the
estimates of the unconditional class probability produced by Eq.~\eqref{eq:general} and Eq.~\eqref{eq:unique}
respectively, indeed can be different. In that example, actually the `total probability' estimate made by means of Eq.~\eqref{eq:general}
is better than the estimate by means of Eq.~\eqref{eq:unique} (but still quite poor) -- although we have argued above that conceptually
the `total odds' is more convincing. Hence, it is not clear whether `total probability' or `total odds' is
better for the estimation of unconditional class probabilities.

However, for an important special case of the probability measure $\mathrm{P}_1$ `total odds' appears to be a
more natural approach to the estimation of the unconditional class probabilities than `total probability'. Under Assumption~\ref{as:LR},
define the probability measure $\mathrm{Q}$ on $(\Omega, \mathcal{H})$ by
\begin{equation}\label{eq:Q}
	\mathrm{Q}(H) = q\,\int_H f_A\,d\mu + (1-q)\,\int_H f_{A^c}\,d\mu, \quad H \in \mathcal{H},
\end{equation}
for a fixed $q \in (0,1)$. By Proposition~\ref{pr:cond.prob} (i), $\mathrm{Q}$ is then absolutely continuous with
respect to $\mathrm{P}_0\bigm|\mathcal{H}$. 

Intuitively, $\mathrm{Q}$ is a modification of $\mathrm{P_0}$ with $p_0= \mathrm{P_0}[A]$
replaced by $q$. But note that $\mathrm{Q}[A]$ is undefined because $A\notin \mathcal{H}$ (otherwise the densities $f_A$ and
$f_{A^c}$ could not be positive $\mu$-almost everywhere). Nonetheless, with this intuition in mind it is natural to favour
such extensions $\mathrm{Q}^\ast$ of $\mathrm{Q}$ to any sub-$\sigma$-field of $\mathcal{A}$ containing $A$ that satisfy
\begin{equation}\label{eq:unbiased}
	\mathrm{Q}^\ast(A) \ = \ q.
\end{equation}
`Total odds' as described in Theorem~\ref{th:main} has this property, and hence provides an unbiased estimator of the unconditional
class probability $q$.

\begin{corollary}\label{co:unbiased}
Let Assumption~\ref{as:LR} hold and define the likelihood ratio $\lambda_0$ as in Proposition~\ref{pr:cond.prob}.
Suppose that $\mathrm{P}_0[\lambda_0 =1] < 1$.
Let $\mathrm{P}_1 = Q$ with $Q$ given by \eqref{eq:Q} for some $0 < q <1$. Then $p_1=q$ is the unique
solution of \eqref{eq:unique} in $(0,1)$ and for $\mathrm{P}_1^\ast$ defined as in Theorem~\ref{th:main} (iii) 
it holds that $\mathrm{P}_1^\ast[A] = q$.
\end{corollary}
See Section~\ref{se:proofs} for the proof of Corollary~\ref{co:unbiased}. In contrast to `total odds', the `total
probability' extension of $\mathrm{Q}$ as given by \eqref{eq:general} does not satisfy \eqref{eq:unbiased} for $q \not= p_0$.
This follows from the next proposition.

\begin{proposition}\label{pr:biased}
Under Assumption~\ref{as:LR}, define the probability measure $\mathrm{Q}$ by \eqref{eq:Q}. Then it holds that
\begin{equation}\label{eq:bias}
|q-p_0| \int \min(f_A, f_{A^c})\,d\mu \ \le \ \bigm| \int \mathrm{P_0}[A\,|\,\mathcal{H}]\,d Q - q\bigm|
\ \le \ |q-p_0|.	
\end{equation}
\end{proposition}
See Section~\ref{se:proofs} for the proof of Proposition~\ref{pr:biased}. The case $f_A = f_{A^c}$ shows that
both inequalities in \eqref{eq:bias} are sharp. As $\int \min(f_A, f_{A^c})\,d\mu$ is
a measure of the classifier's discriminatory power (Bayesian error rate), Proposition~\ref{pr:biased}
suggests that the bias of the estimate of $q$ is the smaller the more powerful the classifier is. 

Interestingly enough, there is a slightly different estimation problem for which the practical performance of 
`total odds' is clearly superior to `total probability'. This problem is the estimation of conditional class probabilities if
targets for the unconditional class probabilities are independently given. \citet[][Chapter~4, Section ``Estimating
the Prior Probabilities'']{Bohn&Stein} describe the problem and two standard solution approaches in
the context of credit rating systems.

Under Assumption~\ref{as:general} the new problem is described as follows:
\begin{itemize}
	\item An estimate (target) $0 < \mathrm{P}^\ast_1[A] < 1$ for an event $A \in \mathcal{A}\backslash\mathcal{H}$ is given.
	Possibly it was produced in a separate, independent estimation exercise. The problem is to construct conditional
	probabilities
	$\mathrm{P}_1^\ast[A\,|\,\mathcal{H}]$ such that 
	\begin{equation}\label{eq:new}
		\mathrm{P}^\ast_1[A] = \mathrm{E_1}\bigl[\mathrm{P}_1^\ast[A\,|\,\mathcal{H}]\bigr].
	\end{equation}
	\item Again, ideally the estimate should be meaningful in the sense of being based on an extension of 
	$\mathrm{P}_1$ to any $\sigma$-field containing $A$, based on observations as given by $(\Omega, \mathcal{A},
	\mathrm{P}_0)$.
\end{itemize}	
The simplest, `total probability' approach to solving Eq.~\eqref{eq:new}	is by setting
\begin{equation}
	\mathrm{P}_1^\ast[A\,|\,\mathcal{H}] = \frac{\mathrm{P}^\ast_1[A]}%
		{\mathrm{E}_1\bigl[\mathrm{P}_0[A\,|\,\mathcal{H}]\bigr]}\,\mathrm{P}_0[A\,|\,\mathcal{H}].
\end{equation}
This approach is unsatisfactory because it is possible that $\mathrm{P}_1^\ast[A\,|\,\mathcal{H}] > 1$
with positive probability under $\mathrm{P}_1$. Of course, this could be interpreted as evidence
of incompatibility as in the case of violation of the likelihood ratio condition in Theorem~\ref{th:main} (i).
\citet{Bohn&Stein} present an alternative approach which uses the `change of base rate' theorem  
\citep[][Theorem~2]{Elkan01}. However, the solution by that approach in general does not solve \eqref{eq:new} 
because in practice often the outcome is $\mathrm{P}^\ast_1[A] \not= \mathrm{E_1}\bigl[\mathrm{P}_1^\ast[A\,|\,\mathcal{H}]\bigr]$.

An alternative estimation approach suggested by \citet[][Section 4.2, ``scaled likelihood ratio'']{Tasche2013}
uses Theorem~\ref{th:main}:
\begin{itemize}
	\item Let $p_1 \stackrel{\text{def}}{=} \mathrm{P}^\ast_1[A]$. Solve then the following equation for $c$:
	\begin{equation}\label{eq:solve}
		1 = \mathrm{E_1}\left[\frac{1}{p_1 + (1-p_1)\,c\,\lambda_0}\right].
	\end{equation}
If $\lambda_0$ is non-constant there is a unique solution $c>0$ of Eq.~\eqref{eq:solve}. 
\item Since $0 < p_1 < 1$, Theorem~\ref{th:main} (i) then implies
$$\frac1{\mathrm{E}_1[\lambda_0]} < c < \mathrm{E}_1\left[\tfrac1{\lambda_0}\right].$$
\item Moreover, if the measure $\mathrm{P}_1^\ast$ is defined with $\lambda_0$ replaced
by $c\,\lambda_0$, Theorem~\ref{th:main} (iii) implies that the solution is meaningful because
it results in a proper extension of 
	$\mathrm{P}_1$ to a $\sigma$-field containing $A$.
\item By Theorem~\ref{th:main} (iv), the resulting estimate of the conditional probability 
$\mathrm{P}_1^\ast[A\,|\,\mathcal{H}]$ is as follows:
\begin{equation}
		\mathrm{P}_1^\ast[A\,|\,\mathcal{H}] = \frac{p_1}{p_1 + (1-p_1)\,c\,\lambda_0}.
	\end{equation}
\end{itemize}
With a view on Eq.~\eqref{eq:gen.odds}, the `scaled likelihood ratio' approach could also be called `total
odds' approach.
Results from an estimation exercise on real-world data presented in \citet{Tasche2013} suggest that
`total odds' in general provides better solutions of problem~\eqref{eq:new} than `total probability'. 

\section{Related work}
\label{sec:RelatedResearch}

\citet{saerens2002adjusting} assumed that the marginal distribution $\mathrm{P}_1$ in Assumption~\ref{as:general}
was given by a mixture distribution like in \eqref{eq:Q}. They suggested estimating the parameter $q$ with a maximum
likelihood approach. To describe their proposal in more detail, suppose there is a sample $\omega_1, \ldots, \omega_n$
of independent observations under $\mathrm{P}_1 = Q$. The likelihood function $L$ is then given by
\begin{equation}\label{eq:likelihood}
	L(q, \omega_1, \ldots, \omega_n) \ = \ \prod_{i=1}^n \bigl(q\,f_A(\omega_i) + (1-q)\,f_{A^c}(\omega_i)\bigr).
\end{equation}
With $\lambda_0 = f_{A^c} / f_A$ as in Proposition~\ref{pr:cond.prob} and Theorem~\ref{th:main}, one then obtains for the
log-likelihood function 
$$
\log\bigl(L(q, \omega_1, \ldots, \omega_n)\bigr)
\ =\ \sum_{i=1}^n \log\bigl(f_A(\omega_i)\bigr) + \sum_{i=1}^n \log\bigl(q + (1-q)\,\lambda_0(\omega_i)\bigr).$$
This implies
$$
\frac{\partial}{\partial\,q} \log\bigl(L(q, \omega_1, \ldots, \omega_n)\bigr) \ = \
\sum_{i=1}^n \frac{1-\lambda_0(\omega_i)}{q + (1-q)\,\lambda_0(\omega_i)}.$$
Equating the derivative to $0$ as a necessary condition for a maximum gives
\begin{align}
0 & = \sum_{i=1}^n \frac{1-\lambda_0(\omega_i)}{q + (1-q)\,\lambda_0(\omega_i)} \notag\\
& = \sum_{i=1}^n \frac{1}{q + (1-q)\,\lambda_0(\omega_i)} -
\frac 1{1-q} \sum_{i=1}^n \frac{(1-q)\,\lambda_0(\omega_i)}{q + (1-q)\,\lambda_0(\omega_i)}\notag\\
& = \sum_{i=1}^n \frac{1}{q + (1-q)\,\lambda_0(\omega_i)} - \frac n{1-q}  +
	\frac q{1-q} \sum_{i=1}^n \frac{1}{q + (1-q)\,\lambda_0(\omega_i)} \notag\\
\iff & 1 \ = \ \frac 1 n \sum_{i=1}^n \frac{1}{q + (1-q)\,\lambda_0(\omega_i)}. \label{eq:empiric}
	\end{align}
Equation \eqref{eq:empiric} is \eqref{eq:unique} with $p_1$ replaced by $q$ and 
$\mathrm{P_1}[H] = \frac 1 n \sum_{i=1}^n \mathbf{1}_H(\omega_i)$, $H\in \mathcal{H}$, the empirical distribution
associated with the sample $\omega_1, \ldots, \omega_n$. This observation shows that the sample version of the
total odds estimator is identical with the maximum likelihood estimator of \citet{saerens2002adjusting}.

Based on Theorem~\ref{th:main}, therefore, we have identified a sufficient and necessary condition for the
maximum likelihood estimator to exist (in the binary classification setting). Moreover, the derivation of
\eqref{eq:empiric} shows that the maximum likelihood estimator works for any model where the ratio of
the conditional class densities equals the relative odds $\lambda_0$. Note that \citet[][Eq.~(9)]{duPlessis2014110}
derived \eqref{eq:empiric} but did not discuss the existence of solutions.


\section{Proofs}
\label{se:proofs}

The proof of Theorem~\ref{th:main} is mainly based on the following lemma that generalises Theorem~3.3 of \citet{Tasche2013}.
\begin{lemma}\label{le:main}
Let $X >0$ be a random variable such that $\mathrm{P}[X=1] < 1$. Then there exists a solution $0 \le p < 1$ to the equation
\begin{equation}\label{eq:lemma}
	\mathrm{E}\left[\frac{1}{p+(1-p)\,X}\right] = 1
\end{equation}
if and only if $\mathrm{E}[X] > 1$ and $\mathrm{E}\bigl[X^{-1}\bigr] \ge 1$. If there is a solution 
$0 \le p < 1$ to Eq.~\eqref{eq:lemma} it is unique. The unique solution is $p=0$ if and only
if $\mathrm{E}\bigl[X^{-1}\bigr] = 1$.
\end{lemma}

\begin{figure}[t!p]
\caption{Illustration for the proof of Lemma~\ref{le:main}. The three
possibilities for the shape of the graph of the function $F$ defined by \eqref{eq:F}.}
\label{fig:1}
\begin{center}
\ifpdf
	\includegraphics[width=10cm]{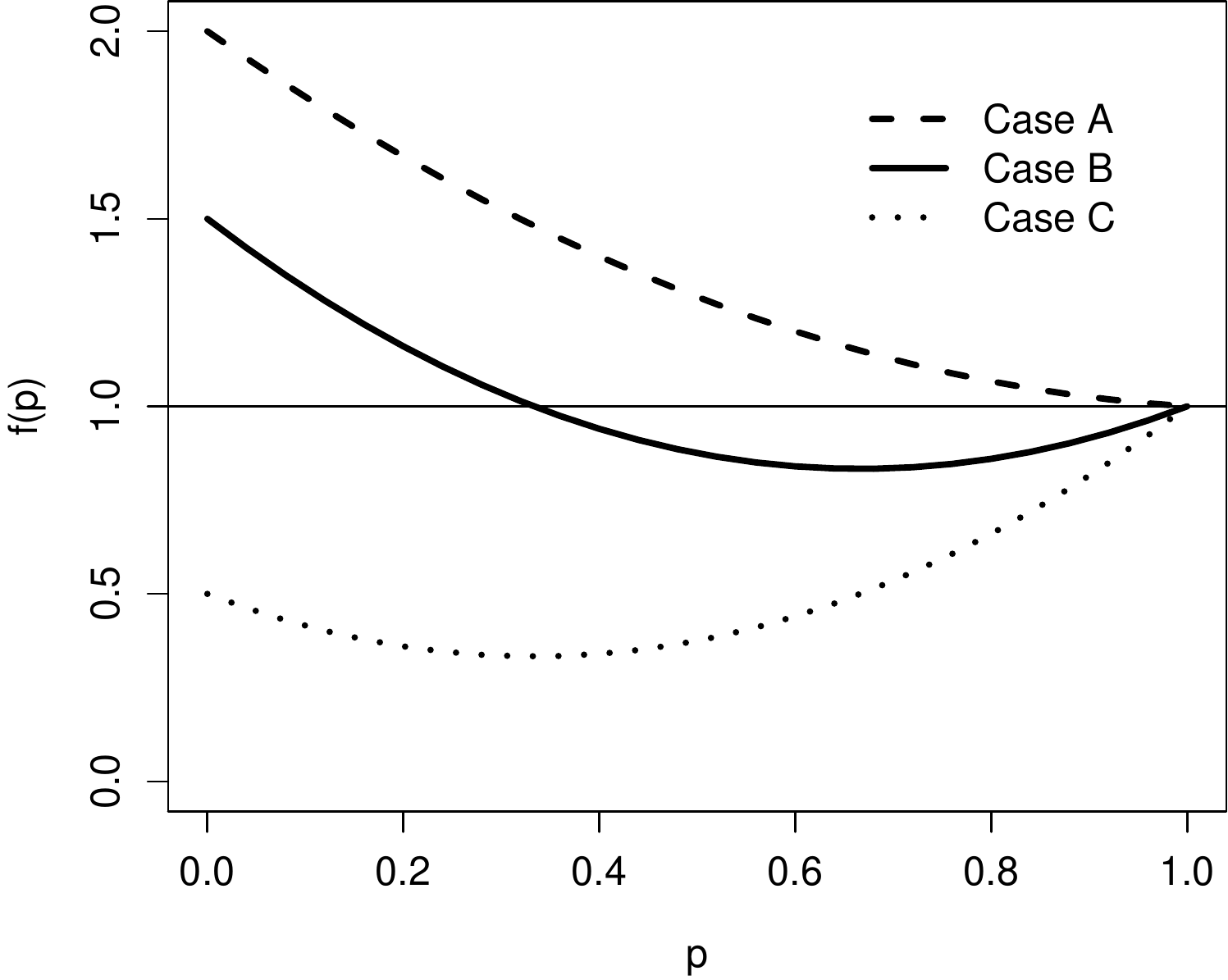}
\fi
\end{center}
\end{figure}

\textbf{Proof.} In principle, the proof in this case is the same as the proof of Theorem~3.3 of \citet{Tasche2013}.
However, we have to take care of the possibility that $\mathrm{E}[X] = \infty$ or $\mathrm{E}\bigl[X^{-1}\bigr] = \infty$.
Define the function $F: [0,1] \to (0,\infty]$, $p \mapsto F(p)$ by
\begin{equation}\label{eq:F}
	F(p) = \mathrm{E}\left[\frac{1}{p+(1-p)\,X}\right].
\end{equation}
Then for $0 < p \le 1$ we have $F(p) \le \frac1p < \infty.$ Solely for $p=0$ it may happen that $F(0) = \infty$, depending on
whether or not $X^{-1}$ is integrable. By the dominated convergence theorem $F(p)$ is continuous in $(0,1]$. If 
$\mathrm{E}[X^{-1}] < \infty$ then again by the dominated convergence theorem $F(p)$ is also continuous in $p=0$ since 
$\frac{1}{p+(1-p)\,X} \le \max(X^{-1},1)$. 
However, Fatou's lemma implies that  $F(p) \xrightarrow{p \rightarrow 0} \mathrm{E}[X^{-1}]$ even if $\mathrm{E}[X^{-1}] = \infty$.

The function $p \mapsto f_X(p) = \frac{1}{p+(1-p)\,X}$ is twice continuously differentiable in $(0,1)$ with 
\begin{align}
	f_X'(p) = \frac{X-1}{(p+(1-p)\,X)^2},\label{eq:deriv}\\
	f_X''(p) = \frac{2\,(X-1)^2}{(p+(1-p)\,X)^3}.\notag
\end{align}
For fixed $p\in (0,1)$ the random variable $f_X'(p)$ is integrable because it holds that
\begin{equation*}
|f_X'(p)| \le \frac{1}{p^2} + \frac 1 p\, \frac{X}{p + (1-p)\,X} =
\frac{1}{p^2} + \frac 1 {p\,(1-p)} \left(\frac{p + (1-p)\,X}{p + (1-p)\,X} - \frac{p}{p + (1-p)\,X}\right)
\le \frac{1}{p^2} + \frac 1 {p\,(1-p)}.
\end{equation*}
Hence it follows from the dominated convergence theorem that also $F$ as defined by \eqref{eq:F} is 
continuously differentiable in $(0,1)$. Moreover, since $f_X''(p) > 0$ on $\{X \not=1\}$ and 
$\mathrm{P}[X=1] < 1$ we obtain that the derivative of $F$ is strictly increasing for $0 < p < 1$ (strict convexity).
Together with the (quasi-)continuity of $F$ this observation implies uniqueness of any solution $0 \le p <1$ to
\eqref{eq:lemma} if there is one.

The strict convexity of $F$ implies that  the graph of $F$ must look like one of the three stylised graphs in Figure~\ref{fig:1}.
Only in case~B is there a solution to Eq.~\eqref{eq:lemma} other than $p=1$. Case~B is characterised by the two conditions
\begin{align*}
\lim_{p\to 0} F(p) & \ge 1 \quad\text{and}\\
\lim_{p\to 1} F'(p) & > 0.
	\end{align*}
We have seen above that $\lim_{p\to 0} F(p) = \mathrm{E}[X^{-1}]$.
Eq.~\eqref{eq:deriv} implies by means of a combination of the dominated convergence
theorem and Fatou's lemma that for both the case $\mathrm{E}[X] < \infty$ and the case $\mathrm{E}[X] = \infty$ we have
$$\lim_{p\to 1} F'(p) = \lim_{p\to 1} \mathrm{E}[f'(p)] = \mathrm{E}[X] - 1.$$
This proves the existence part of the lemma. The criterion for the solution to \eqref{eq:lemma} to be
$p=0$ also follows from $\lim_{p\to 0} F(p) = \mathrm{E}[X^{-1}]$. \hfill \qed

\textbf{Proof of Theorem~\ref{th:main}.} (i) is an immediate conclusion from Lemma~\ref{le:main}. Since 
$\bigl\{(A\cap H) \cup (A^c\cap G): H, G \in \mathcal{H}\bigr\}$ is a $\sigma$-field (ii) follows from
the observation
$$\mathcal{H} \cup \{A\} \subset \bigl\{(A\cap H) \cup (A^c\cap G): H, G \in \mathcal{H}\bigr\} \subset 
\sigma(\mathcal{H} \cup \{A\}).$$
We begin the proof of (iii) with another lemma.
\begin{lemma}\label{le:null}
Let $H\in\mathcal{H}$. Then
\begin{align*}
A \cap H = \emptyset &\ \Rightarrow \mathrm{E_1}\left[\mathbf{1}_{H}\,\frac{p_1}{p_1 + (1-p_1)\,\lambda_0}\right] = 0,\\
A^c \cap H = \emptyset &\ \Rightarrow	\mathrm{E_1}\left[\mathbf{1}_{H}\,\frac{(1-p_1)\,\lambda_0}{p_1 + (1-p_1)\,\lambda_0}\right] = 0.
	\end{align*}
\end{lemma}
\textbf{Proof of Lemma~\ref{le:null}.} Denote by $\varphi$ 
any $\mathcal{H}$-measurable density of $\mathrm{P}_1$ with respect to $\mathrm{P}_0$. Proposition~\ref{pr:cond.prob} (ii) then implies
\begin{align*}
	\mathrm{E_1}\left[\mathbf{1}_{H}\,\frac{p_1}{p_1 + (1-p_1)\,\lambda_0}\right] & =
\mathrm{E_0}\left[\varphi\,\mathbf{1}_{H}\,\frac{\frac{p_1}{p_0}\,\mathrm{P}_0[A\,|\,\mathcal{H}]}
					{\frac{p_1}{p_0}\,\mathrm{P}_0[A\,|\,\mathcal{H}] + 
					\frac{1-p_1}{1-p_0}\,\mathrm{P}_0[A^c\,|\,\mathcal{H}]}\right]\\
			& = \frac{p_1}{p_0}\,\mathrm{E_0}\left[\varphi\,\mathbf{1}_{H\cap A}\,\frac{1}
					{\frac{p_1}{p_0}\,\mathrm{P}_0[A\,|\,\mathcal{H}] + 
					\frac{1-p_1}{1-p_0}\,\mathrm{P}_0[A^c\,|\,\mathcal{H}]}\right]\\
					& = 0.
\end{align*}	
The proof of the second implication in Lemma~\ref{le:null} is almost identical. \hfill \qed

\textbf{Proof of Theorem~\ref{th:main} continued.} Let $B \in \mathcal{H}^A$ with
$$B = (A\cap H_1) \cup (A^c\cap G_1) = (A\cap H_2) \cup (A^c\cap G_2),$$
for some $H_1, H_2, G_1, G_2 \in \mathcal{H}$. Then it follows that
	\begin{gather*}
A\cap H_1 = A\cap H_2 = A\cap H_1 \cap H_2 \ \text{and}\ 
A^c\cap G_1 = A^c\cap G_2 = A^c\cap G_1 \cap G_2.
\end{gather*}
Hence $A\cap (H_1\backslash H_2) = \emptyset = A\cap (H_2\backslash H_1)$ and
$A^c\cap (G_1\backslash G_2) = \emptyset = A^c\cap (G_2\backslash G_1)$. 
Lemma~\ref{le:null} now implies that $\mathrm{P}_1^\ast$ is well-defined because it holds for any sets $M_1$, $M_2$ that
$$\mathbf{1}_{M_1} = \mathbf{1}_{M_1\cap M_2} +  \mathbf{1}_{M_1\backslash M_2} 
\ \text{and}\ \mathbf{1}_{M_2} = \mathbf{1}_{M_1\cap M_2} +  \mathbf{1}_{M_2\backslash M_1}.$$

The properties $\mathrm{P}_1^\ast[\emptyset] =0$, $\mathrm{P}_1^\ast[\Omega] =1$ and 
$\mathrm{P}_1^\ast[H] = \mathrm{P}_1[H]$ for $H \in \mathcal{H}$ are obvious. 
Finite addivity of $\mathrm{P}_1^\ast$ follows from Lemma~\ref{le:null} because
$$B_i = (A \cap H_i) \cup (A^c \cap G_i),\ i=1,2\ \text{with}\ B_1 \cap B_2 = \emptyset$$
implies $A \cap H_1\cap H_2 = \emptyset =  A^c \cap G_1 \cap G_2$ and
$$B_1 \cup B_2 = \bigl(A \cap (H_1\cup H_2)\bigr) \cup \bigl(A^c \cap (G_1 \cup G_2)\bigr).$$
To complete the
proof of (iii) we have to show that $\mathrm{P}_1^\ast$ is $\sigma$-continuous in $\emptyset$, i.e.
\begin{equation}\label{eq:conv}
	\lim_{n\to\infty} \mathrm{P}_1^\ast[B_n] = 0,
\end{equation}
for any $B_1 \supset B_2 \supset \ldots$ with $\bigcap_{n=1}^\infty B_n = \emptyset$. Let
$(B_n)$ be such a sequence in $\mathcal{H}^A$ with representation
$B_n =(A \cap H_{n}) \cup (A^c \cap G_{n})$,
for sequences $(H_{n})$, $(G_{n})$ in $\mathcal{H}$. 
Note that 
$$\frac{p_1}{p_0}\,\mathrm{P}_0[A\,|\,\mathcal{H}] + 
					\frac{1-p_1}{1-p_0}\,\mathrm{P}_0[A^c\,|\,\mathcal{H}] \ge
					\min\left(\frac{p_1}{p_0}, \frac{1-p_1}{1-p_0}\right).$$
					Therefore, similarly to
the proof of Lemma~\ref{le:null} we see that
\begin{align*}
\mathrm{P}_1^\ast[B_n] & \le \frac{\max\bigl(p_1/p_0, (1-p_1)/(1-p_0)\bigr)}
{\min\bigl(p_1/p_0, (1-p_1)/(1-p_0)\bigr)}\,
\mathrm{E}_0\left[\varphi\,\bigl(\mathbf{1}_{H_{n}}\,\mathrm{P}_0[A\,|\,\mathcal{H}] +
\mathbf{1}_{G_{n}}\,\mathrm{P}_0[A^c\,|\,\mathcal{H}]\bigr)\right]\\
& = \frac{\max\bigl(p_1/p_0, (1-p_1)/(1-p_0)\bigr)}
{\min\bigl(p_1/p_0, (1-p_1)/(1-p_0)\bigr)}\, \mathrm{E}_0[\varphi\,\mathbf{1}_{B_n}],
	\end{align*}
where $\varphi$ is an $\mathcal{H}$-measurable density as in Lemma~\ref{le:null}	.
By the dominated convergence theorem, Eq.~\eqref{eq:conv} follows. 

With regard to (iv), observe that by the definition of $\mathrm{P}^\ast_1$ and the fact 
that $\mathrm{P}^\ast_1\bigm|\mathcal{H} = \mathrm{P}_1$ it holds for $H \in \mathcal{H}$ that
$$\mathrm{E_1}\left[\mathbf{1}_{H}\,\frac{p_1}{p_1 + (1-p_1)\,\lambda_0}\right] = 
\mathrm{P}^\ast_1[A\cap H] = \mathrm{E}_1^\ast\bigl[\mathbf{1}_H\,\mathrm{P}_1^\ast[A\,|\,\mathcal{H}]\bigr]
= \mathrm{E}_1\bigl[\mathbf{1}_H\,\mathrm{P}_1^\ast[A\,|\,\mathcal{H}]\bigr].$$
This implies (iv) because $\lambda_0$ is $\mathcal{H}$-measurable.
\hfill \qed

\textbf{Proof of Corollary~\ref{co:unique}.} Note that \eqref{eq:equiv} is equivalent to 
\begin{align*}
\widetilde{\mathrm{P}}_1[A\,|\,\mathcal{H}] & =  \frac{\widetilde{\mathrm{P}}_1[A]}{\widetilde{\mathrm{P}}_1[A] +
(1-\widetilde{\mathrm{P}}_1[A])\,\lambda_0.}\\
\intertext{This implies}
1 & = \mathrm{E_1}\left[\frac{1}{\widetilde{\mathrm{P}}_1[A] +
(1-\widetilde{\mathrm{P}}_1[A])\,\lambda_0}\right].
	\end{align*}
Therefore, by Theorem~\ref{th:main} (i) we can conclude that
$\widetilde{\mathrm{P}}_1[A] = p_1$. By Theorem~\ref{th:main} (iv), it follows that 
$\widetilde{\mathrm{P}}_1[A\,|\,\mathcal{H}] = \mathrm{P}^\ast_1[A\,|\,\mathcal{H}]$ and hence
$$\widetilde{\mathrm{P}}_1[A \cap H] \ = \ \mathrm{E}_1\bigl[\mathbf{1}_H\,\mathrm{P}^\ast_1[A\,|\,\mathcal{H}]\bigr] \ =\ 
\mathrm{P}^\ast_1[A \cap H], \quad H \in \mathcal{H}.$$
This implies $\widetilde{\mathrm{P}}_1 = \mathrm{P}_1^\ast$ because $A \cap \mathcal{H}$ is a $\cap$-stable generator of
$\mathcal{H}^A$. \hfill \qed

\textbf{Proof of Corollary~\ref{co:unbiased}.} Observe that $q\,f_A + (1-q)\,f_{A^c}$ is a $\mu$-density of $\mathrm{P}_1 = Q$.
This implies
$$\mathrm{E}_1\bigl[\frac1{p_1+(1-p_1)\,\lambda_0}\bigr] =
\int \frac{q\,f_A + (1-q)\,f_{A^c}}{p_1\,f_A + (1-p_1)\,f_{A^c}}\,f_A\,d\mu = 1,$$
if we choose $p_1 = q$. As $f_A$ and $f_{A^c}$ are positive $\mu$-almost everywhere, $\mathrm{P}_0[\lambda_0 =1] < 1$ implies
$\mathrm{P}_1[\lambda_0 =1] = \mathrm{Q}[\lambda_0 =1]< 1$. 
By Theorem~\ref{th:main} (i), hence the moment conditions on $\lambda_0$ are satisfied and 
there is no other solution to \eqref{eq:unique} than $q$. From this it follows that $\mathrm{P}_1$ can be extended to
$\mathcal{H}^A$ as defined in Theorem~\ref{th:main} (ii) and that the extension satisfies $\mathrm{P}_1^\ast[A] = q$.
\hfill \qed

\textbf{Proof of Proposition~\ref{pr:biased}.} If $\mu(f_A \not= f_{A^c}) = 0$ then
all three parts of \eqref{eq:bias} equal $|q-p_0|$. Suppose now that $\mu(f_A \not= f_{A^c}) > 0$.
By Proposition~\ref{pr:cond.prob} (ii) we can calculate as follows:
\begin{align}
\int \mathrm{P_0}[A\,|\,\mathcal{H}]\,d Q - q & = p_0 \int f_A\,\frac{q\,f_A + (1-q)\,f_{A^c}}{p_0\,f_A + (1-p_0)\,f_{A^c}}\,d\mu - q\notag\\
& = \int f_A\,\frac{p_0\,q\,f_A + p_0\,(1-q)\,f_{A^c} - q\,p_0\,f_A -q\, (1-p_0)\,f_{A^c}}{p_0\,f_A + (1-p_0)\,f_{A^c}}\,d\mu\notag\\
& = (p_0 - q) \int \frac{f_A\,f_{A^c}}{p_0\,f_A + (1-p_0)\,f_{A^c}}\,d\mu. \label{eq:estimate}
	\end{align}
Observing that $\frac{f_A\,f_{A^c}}{p_0\,f_A + (1-p_0)\,f_{A^c}} \ge \min(f_A, f_{A^c})$ we obtain the first 
inequality in \eqref{eq:bias}.  With regard to the second inequality, define a probability measure 
$\mathrm{P}$ on $(\Omega, \mathcal{H})$ by 
$$\mathrm{P}[H] \ =\ \int_H f_{A^c}\,d\mu, \quad H \in \mathcal{H}.$$
With $X = \frac{f_{A^c}}{f_A}$ then if follows for all $p \in [0,1]$ that
\begin{equation}\label{eq:apply}
\int \frac{f_A\,f_{A^c}}{p\,f_A + (1-p)\,f_{A^c}}\,d\mu \ =\ \mathrm{E}\left[\frac1{p + (1-p)\,X}\right].
\end{equation}
 Note that 
$\mathrm{E}\left[\tfrac1{X}\right]\ = \ \int f_A \,d\mu = 1$.
In addition, since $f_{A}$ is positive $\mu(f_A \not= f_{A^c}) > 0$ implies 
$\mathrm{P}[X=1] < 1$. Hence we can apply Lemma~\ref{le:main} to conclude that $p=0$ is the only $p\in[0,1)$
such that $\mathrm{E}\left[\frac1{p + (1-p)\,X}\right] = 1$. The proof of Lemma~\ref{le:main} shows that in
this case for $0 < p < 1$ we have
$$1\ >\ \mathrm{E}\left[\frac1{p + (1-p)\,X}\right].$$
By \eqref{eq:apply} and \eqref{eq:estimate}, the second inequality in \eqref{eq:bias} follows.
\hfill \qed

Note that \eqref{eq:estimate} could be rearranged in order to construct an unbiased estimator of $q$.



\addcontentsline{toc}{section}{References}

\end{document}